\documentclass[11pt]{article}
\usepackage[letterpaper,hmargin=1in,vmargin=1.25in]{geometry}
\usepackage{latexsym, amsfonts, amssymb, amsmath, amsthm, graphicx,hyperref}

\newcommand{\charf}[1]{\mbox{\raise.48ex\hbox{$\chi$}$_{#1}$}}
\def\~{\tilde }
\def\C{{\cal C}}

\def\R{{\mathbb R}}
\def\ZZ{{\mathbb Z}}
\def\Q{{\sf Q}}
\def\P{{\sf P}}
\def\p{{\sf p}}
\def\I{{\sf I}}
\def\E{{\sf E}}

\def\1{{\sf 1}}
\def\K{{\sf K}}

\theoremstyle{plain}
\newtheorem{theorem}{Theorem}[section]
\newtheorem{lemma}[theorem]{Lemma}
\newtheorem{corollary}[theorem]{Corollary}

\newtheorem{remark}[theorem]{Remark}

\theoremstyle{definition}
\newtheorem{definition}[theorem]{Definition}
\newtheorem{example}[theorem]{Example}

\begin{document}

\title{Sharp edge, vertex, and mixed Cheeger type inequalities\\ for finite Markov kernels}

\author {Ravi Montenegro \thanks{Department of Mathematical Sciences, University of
      Massachusetts Lowell, Lowell, MA 01854, ravi\_montenegro@uml.edu; supported in part by
      a VIGRE grant at the Georgia Institute of Technology.}
}

\date{}

\maketitle
 
\begin{abstract}
\noindent
We show how the evolving set methodology of Morris and Peres can be used to show Cheeger inequalities for bounding the spectral gap of a finite Markov kernel. This leads to sharp versions of several previous Cheeger inequalities, including ones involving edge-expansion, vertex-expansion, and mixtures of both. A bound on the smallest eigenvalue also follows.

\vspace{1ex}
\noindent {\bf Keywords} : Markov chain, evolving sets, Cheeger inequality, eigenvalues.
\end{abstract}

\noindent
\section{Introduction}

Given a finite, irreducible reversible Markov kernel $\P$ the Perron-Frobenius theorem guarantees that the matrix $\P$ has a real valued eigenbasis with eigenvalues $1=\lambda_0(\P)\geq\lambda_1(\P)\geq\cdots\geq\lambda_{n-1}(\P)\geq -1$. The spectral gap $\lambda=1-\lambda_1(\P)$ between the largest and second largest eigenvalues, or in the non-reversible case the gap \mbox{$\lambda=1-\lambda_1\left(\frac{\P+\P^*}{2}\right)$} of the additive symmetrization, governs key properties of the Markov chain. Alon \cite{Alon86.1}, Lawler and Sokal \cite{LS88.1}, and Jerrum and Sinclair \cite{JS88.1} showed lower bounds on the spectral gap in terms of geometric quantities on the underlying state space $V$, known as Cheeger inequalities. Similarly, in the reversible case Diaconis and Stroock \cite{DS91.1} showed a lower bound on $1+\lambda_{n-1}$, known as Poincar\'e inequalities, which also have a geometric flavor.

Inequalities of both types have played an important role in the study of the mixing times of Markov chains. Conversely, the authors of \cite{MP05.1} used their Evolving set bounds on mixing times to show a Cheeger inequality, although this was later removed as it was weaker than previously known bounds. We improve on their idea and find that our resulting Theorem \ref{thm:cheeger} can be used to show sharp Cheeger-like lower bounds on $\lambda$ and $1+\lambda_{n-1}$, both in the edge-expansion sense of Jerrum and Sinclair, the vertex-expansion notion of Alon, and a mixture of both. The bounds on $\lambda$ typically improve on previous bounds by a factor of two, which is essentially all that can be hoped for as most of our bounds are sharp; the notion of edge-expansion used in our Cheeger inequality for $1+\lambda_{n-1}$ is entirely new.

The paper is organized as follows. In the preliminaries we review some mixing time and Evolving set results. This is followed in Section \ref{sec:general-cheeger} by our main result, an Evolving set generalization of Cheeger's inequality. In Section \ref{sec:cheeger} this is used to show a sharp version of the edge expansion Cheeger Inequality and to improve on vertex-expansion bounds of Alon and of Stoyanov. Similar bounds on $\lambda_{n-1}$, and more generally the second largest magnitude eigenvalue, are found in Section \ref{sec:small_eigenvalue}.

%*********************** Preliminaries **************

\section{Preliminaries} 

Consider a finite ergodic Markov kernel $\P$ (i.e. transition probability matrix) on state space $V$ with stationary distribution $\pi$. This is called lazy if $\P(x,x)\geq 1/2$ for every $x\in V$, and is reversible if $\P^*=\P$ where the time-reversal $\P^*(x,y)=\frac{\pi(y)\P(y,x)}{\pi(x)}$. The ergodic flow from $A\subset V$ to $B\subset V$ is $\Q(A,B)=\sum_{x\in A,y\in B} \pi(x)\P(x,y)$. The total variation distance between distributions  $\sigma$ and $\pi$ is $\|\sigma-\pi\|_{TV}=\frac 12\sum_{x\in V}|\sigma(x)-\pi(x)|$. The rate of convergence of a reversible walk is related to spectral gap \cite{JS88.1,DS91.1} by 
\begin{equation} \label{eqn:mixing-lower}
\frac 12\,(1-\lambda)^n \leq \frac 12\,\lambda_{max}^n \leq \max_{x\in V} \|\p_x^n-\pi\|_{TV}
\leq \frac 12\,\frac{\lambda_{max}^n}{\min_{y\in V}\pi(y)}\,
\end{equation}
where $\p_x^n(y)=\P^n(x,y)$ and $\lambda_{max}=\max\{\lambda_1(\P),|\lambda_{n-1}(\P)|\}$.

Morris and Peres \cite{MP05.1} introduced a new tool for studying the rate of convergence:

\begin{definition}
Given set $A\subset V$ a step of the {\em evolving set process} is given by choosing $u\in[0,1]$ uniformly at random, and transitioning to the set 
$$
A_u = \{y\in V : \Q(A,y) \geq u\,\pi(y)\}=\{y\in V: \P^*(y,A)\geq u\}\,.
$$
The walk is denoted by $S_0$, $S_1$, $S_2$, $\ldots$, $S_n$, with transition kernel $\K^n(A,S)=Prob(S_n=S|S_0=A)$, and expectation $\E_n f(S_n)=\sum_{S_n\subset V} K^n(S_0,S_n)\,f(S_n)$.
\end{definition}

The main result of \cite{MP05.1} is a bound on the rate of convergence in terms of Evolving sets:

\begin{lemma} \label{lem:MP}
If $S_0=\{x\}$ for some $x\in V$ then
$$
\|\p_x^n-\pi\|_{TV} \leq \frac{1}{2\,\pi(x)}\E_n\sqrt{\min\{\pi(S_n),1-\pi(S_n)\}}\,.
$$
\end{lemma}

A few easy lemmas of theirs will be required for our work, both of which the interested reader should have little trouble in showing. First, a Martingale relation:

\begin{lemma} \label{lem:martingale} 
If $A\subset V$ then
$$
\int_0^1 \pi(A_u)\,du=\pi(A)\,.
$$
\end{lemma}

Note from the definition that for a lazy walk $A_u\subset A$ if $u>1/2$, while $A_u\supset A$ if $u\leq 1/2$. The gaps between $A$ and $A_u$ are actually related to ergodic flow:

\begin{lemma} \label{lem:ergodic-flow}
Given a lazy Markov chain, if $A\subset V$ then
$$
\Q(A,A^c)=\int_0^{1/2}(\pi(A_u)-\pi(A))\,du=\int_{1/2}^1(\pi(A)-\pi(A_u))\,du\,.
$$
\end{lemma}

%********************** Generalized Cheeger *********************

\section{A Generalized Cheeger Inequality} \label{sec:general-cheeger}

Recall that a Cheeger inequality is used to bound eigenvalues of the transition kernel in terms of some geometric quantity. 
``The Cheeger Inequality'' in the finite Markov setting generally refers to the bound
\begin{equation} \label{eqn:normal-cheeger}
\lambda \geq 1-\sqrt{1-h^2}\geq \frac{h^2}{2}
\quad\textrm{where}\quad
h = \min_{\substack{A\subset V,\\ \pi(A)\leq 1/2}} \frac{\Q(A,A^c)}{\pi(A)}\,.
\end{equation}
The quantity $h$ is known as the Cheeger constant, or Conductance, and measures how quickly the walk expands from a set.
We now show a generalization of the Cheeger inequality which is expressed in terms of Evolving sets. The Cheeger constant is replaced by $f$-congestion:

\begin{definition}
If $f:[0,1]\to\R_+$ and $A\subset V$ the $f$-congestion of $A$ is given by
$$
\C_f(A)=\frac{\int_0^1 f(\pi(A_u))\,du}{f(\pi(A))}\,.
$$
The $f$-congestion is given by $\C_f = \max_{A\subset V} \C_f(A)$.
\end{definition}

Small $f$-congestion corresponds to a rapid change in set size of the Evolving set process. In \cite{MT06.1} it is found that to study many measures of convergence rate (total variation, relative entropy, chi-square, etc.) there correspond appropriate choice of $\C_f$. The $f$-congestion is thus closely related to convergence of Markov chains, which in part explains why our main result holds:

\begin{theorem} \label{thm:cheeger}
Given a finite, irreducible, reversible Markov chain, and $f:[0,1]\to\R_+$ then
$$
\lambda \geq 1-\lambda_{max} \geq 1-\C_f\,.
$$
If $\forall a\in(0,1/2):f(a)\leq f(1-a)$ then it suffices to let $\C_f = \max_{\pi(A)\leq 1/2} \C_f(A)$.
\end{theorem}

\begin{proof}
Given $x\in V$ let $S_0=\{x\}$ and \mbox{$M=\max_{\pi(A)\neq 0,\,1} \frac{\sqrt{\pi(A)}}{2\,g(\pi(A))}$}, where $g:\,[0,1]\to\R_+$ is some function to be defined later. Then by Lemma \ref{lem:MP},
\begin{eqnarray*}
\|\p_x^n-\pi\|_{TV} 
%  &\leq&  \frac{1}{\pi(x)}\,\E_n \sqrt{\pi(S_n^{\#})} \\
   &\leq& \frac{M}{\pi(x)}\,\E_n\,g(\pi(S_n)) \\
   &\leq& \frac{M}{\pi(x)}\,\E_{n-1}\,g(\pi(S_{n-1}))\,\C_g(S_{n-1}) \\
   &\leq& \frac{M\,g(\pi(x))}{\pi(x)}\,\C_g^n\,.
\end{eqnarray*}
The final inequality followed from $\C_g(S_{n-1})\leq\C_g$, induction, and $S_0=\{x\}$. 

But then, by equation \eqref{eqn:mixing-lower},
$$
\lambda_{max} \leq \sqrt[n]{2 \max_x \|\p_x^n-\pi\|_{TV} } 
  \leq \C_g \,\sqrt[n]{2\left(\max_x \frac{M\,g(\pi(x))}{\pi(x)}\right)} \xrightarrow{n\to\infty} \C_g\,.
$$

If $\forall a\in(0,1/2):f(a)\leq f(1-a)$ then let $g(a)=f(\min\{a,1-a\})$, noting that since $\pi((A^c)_u)=\lim_{\delta\to 0} 1-\pi(A_{1-u+\delta})$ then $\int_0^1 g(\pi((A^c)_u))\,du = \int_0^1 g(\pi(A_u))\,du$. Otherwise let $g(a)=f(a)$.
\end{proof}

\begin{remark}
For a non-reversible walk Theorem \ref{thm:cheeger} holds with $1-\lambda_{max}$ replaced by $1-\lambda_*$, where $\lambda_*=\max_{i>0} |\lambda_i|$ is the second largest magnitude (complex-valued) eigenvalue of $\P$. This follows from the related lower bound
$$
\frac 12\,\lambda_*^n \leq \max_{x\in V}\|\p_x^n-\pi\|_{TV}
$$
(see e.g. \cite{MT06.1}). While intriguing, it is unclear if lower bounds on $1-\lambda_*$ have any practical application.
\end{remark}

\begin{remark}
An anonymous reader notes that rather than using lower bounds on variation distance we could instead use the well known-relation $\rho(A)=\lim_{k\to\infty}\|A^k\|^{1/k}$ for spectral radius in terms of a consistent matrix norm satisfying $\|Av\|\leq C\,\|A\|\,\|v\|$. In this case take $A=\P-E$ where $E$ is the matrix with rows all equal to $\pi$, and total variation norm has $C=2$. 
\end{remark}

%**************** Cheeger Inequalities ***************

\section{Cheeger Inequalities} \label{sec:cheeger}

Special cases of Theorem \ref{thm:cheeger} 
%the Generalized Cheeger inequality 
include bounds of the vertex type as in Alon \cite{Alon86.1}, the edge type as in Jerrum and Sinclair \cite{JS88.1}, and mixtures of both. The key to the reduction will be the following lemma:

\begin{lemma} \label{lem:worst-case}
Given a concave function $f:[0,1]\rightarrow \R$ and
two non-increasing functions $g,\,\hat{g}:\,[0,1]\rightarrow[0,1]$
such that $\int_0^1 g(u)\,du=\int_0^1 \hat{g}(u)\,du$ and
$\forall t\in[0,1]:\,\int_0^t g(u)\,du \geq \int_0^t \hat{g}(u)\,du$, 
then
$$
\int_0^1 f\circ g(u)\,du \leq \int_0^1 f\circ \hat{g}(u)\,du\,.
$$
\end{lemma}

\begin{proof}
If $x\geq y,\,\delta> 0$ then $\lambda = 1 - \frac{\delta}{x-y+2\delta}\in[\mbox{$\frac 12$},1]$ with $x=(1-\lambda)\,(y-\delta)+\lambda\,(x+\delta)$ and $y=\lambda\,(y-\delta)+(1-\lambda)\,(x+\delta)$. Concavity of $f$ implies that $f(x)\geq (1-\lambda)\,f(y-\delta)+\lambda\,f(x+\delta)$ and $f(y)\geq \lambda\,f(y-\delta)+(1-\lambda)\,f(x+\delta)$. It follows that
\begin{equation} \label{eqn:bigup_smalldown}
\forall x\geq y,\,\delta\geq 0:\,f(x)+f(y) \geq f(x+\delta) + f(y-\delta)\,.
\end{equation}

The inequality \eqref{eqn:bigup_smalldown} shows that if a bigger value ($x$) is increased
by some $\delta$, while a smaller value ($y$) is decreased by $\delta$, then the sum
$f(x)+f(y)$ decreases.
In our setting, the condition that $\forall t\in[0,1]:\,\int_0^t g(u)\,du \geq \int_0^t \hat{g}(u)\,du$
shows that changing from $\hat{g}$ to $g$ increased the already large values of $\hat{g}(u)$ when $u$ is small,
while the equality $\int_0^1 g(u)\,du=\int_0^1 \hat{g}(u)\,du$ assures that this is canceled out by
an equal decrease in the already small values when $u$ is big. The lemma then follows from \eqref{eqn:bigup_smalldown}.
\end{proof}

It remains to relate the $f$-congestion to the edge or vertex notions of Cheeger constant, and then choose the optimal function $f$. 

%****************** Edge Expansion ********************

\subsection{Edge expansion}

We first consider edge-expansion, i.e. ergodic flow, and in particular derive a bound in terms of the symmetrized Cheeger constant
$$
\~h = \min_{A\subset V} \~h(A)
\quad\textrm{where}\quad
\~h(A) = \frac{\Q(A,A^c)}{\pi(A)\pi(A^c)}\,.
$$
To do this a somewhat stronger bound will be shown, and then a few special cases will be considered, including that of $\~h$.

\begin{corollary} \label{cor:diffiQ}
Given function $f:\,(0,1)\rightarrow\R^+$ such that $f$ and $f''$ are concave, and $\forall a\in(0,1/2):\,f(a)\leq f(1-a)$ then the spectral gap of a finite, irreducible Markov chain satisfies
$$
\lambda \geq \min_{\pi(A)\leq 1/2} \frac{\Q^2(A,A^c)}{-f(\pi(A))/f''(\pi(A))}\,.
$$
Without condition $f(a)\leq f(1-a)$ the result holds with minimum taken over all proper subsets of $V$.
\end{corollary}

\begin{proof}
First consider the reversible, lazy case. By Lemma \ref{lem:ergodic-flow} and the remarks before it, $\Q(A,A^c)$ is the area below $\pi(A_u)$ and above $\pi(A)$, and also above $\pi(A_u)$ and below $\pi(A)$. By Lemma \ref{lem:worst-case} the value $\C_f(A)$ is maximized when $\pi(A_u)=m(u)$ where $m(u)$ is as in the first diagram of Figure \ref{fig:extreme}.

\begin{figure}[ht]
\begin{center}
\includegraphics[height=1.25in]{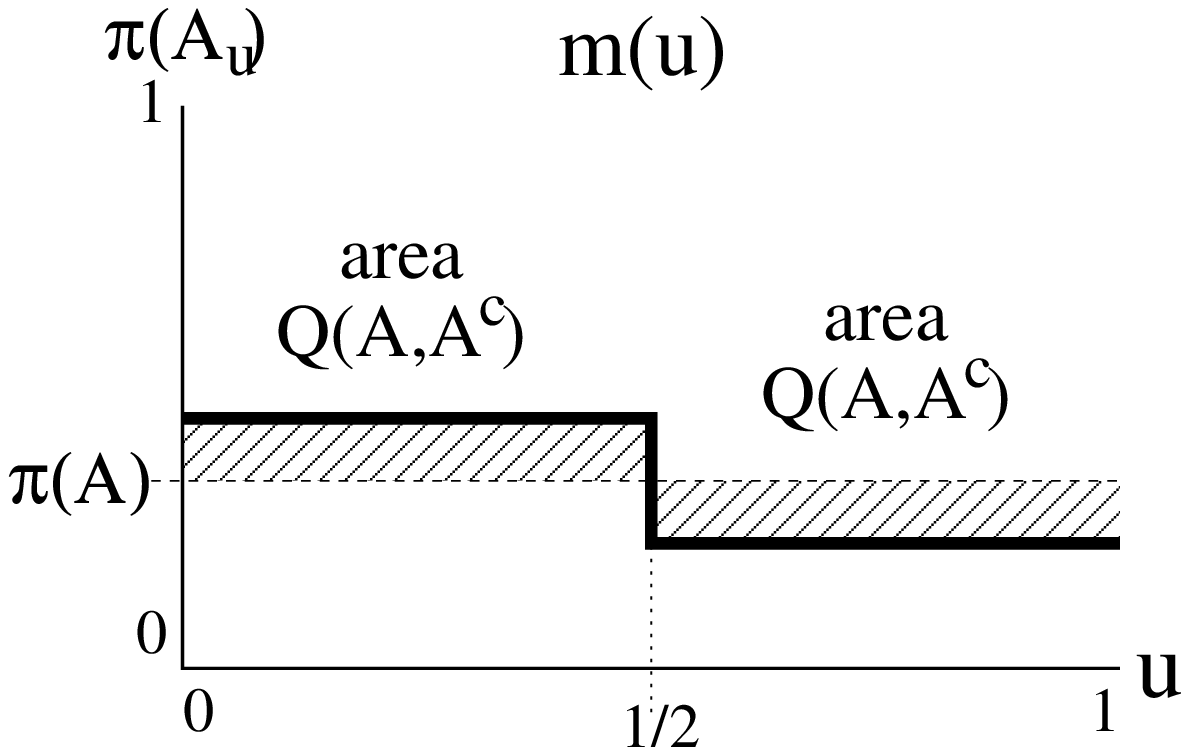}
\hspace{0.5in}
\includegraphics[height=1.25in]{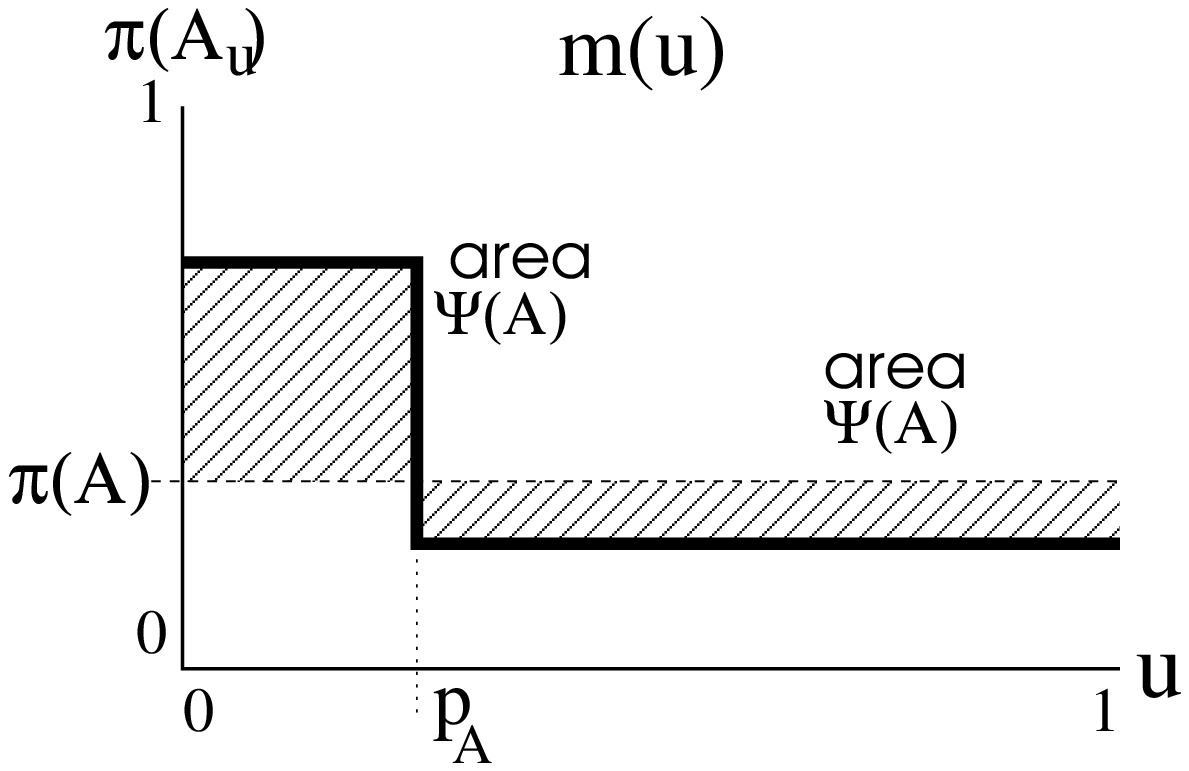}
\caption{Extreme cases,
   $\forall t\in[0,1]:\,\int_0^t \pi(A_u)\,du \geq \int_0^t m(u)\,du$ for lazy and non-lazy walks.} \label{fig:extreme}
\end{center}
\end{figure}

Then
\begin{equation} \label{eqn:lazy-strong-Cheeger}
\lambda \geq 1-\C_f \geq 1-\max_{A\subset V} \frac{f(\pi(A)+2\Q(A,A^c))+f(\pi(A)-2\Q(A,A^c))}{2\,f(\pi(A))}\,.
\end{equation}

In the general case consider the lazy, reversible Markov chain
$\P'=\frac 12\,\left(I + \frac{\P+\P^*}{2}\right)$. Then apply equation \eqref{eqn:lazy-strong-Cheeger} to $\P'$ and observe that $\lambda=2\lambda_{\P'}$ and $\Q(A,A^c)=2\Q_{\P'}(A,A^c)$, to derive the relation
\begin{eqnarray} 
\lambda &\geq& 2\left(1-\max_{A\subset V} \frac{f(\pi(A)+\Q(A,A^c))+f(\pi(A)-\Q(A,A^c))}{2\,f(\pi(A))}\right) \label{eqn:strong-Cheeger} \\
        &\geq& \min_{A\subset V} \frac{\Q(A,A^c)^2}{-f(\pi(A))/f''(\pi(A))}\,. \nonumber
\end{eqnarray}
The second inequality required that $f(x+\delta)+f(x-\delta) \leq 2\,f(x) + f''(x)\,\delta^2$. To show this let $F(y,\delta)=f''(y)\,\delta^2 + 2\,f(y)- f(y+\delta)-f(y-\delta)$. Observe that
$\frac{d}{d\delta} F(y,\delta) = 2\,f''(y)\delta - f'(y+\delta)+f'(y-\delta)=0$ at $\delta=0$,
and $\frac{d^2}{d\delta^2} F(y,\delta) = 2\,f''(y) - (f''(y+\delta)+f''(y-\delta))\geq 0$ as
$f''$ is concave. Hence $F(y,\delta)$ is convex in $\delta$
with minimum at $\delta=0$, and therefore $F(y,\delta)\geq F(y,0)=0$.
\end{proof}
Note that \eqref{eqn:strong-Cheeger} holds even if $f''$ is not concave.

Given $\~h$ then $\Q(A,A^c)\geq\~h \pi(A)\pi(A^c)\ \forall A\subset V$. To apply Corollary \ref{cor:diffiQ} to give a lower bound $\lambda\geq c\~h^2$ we look for a solution to $-f(x)/f''(x)\leq c^{-1}\,x^2(1-x)^2$ where $f,f''$ concave and $c\in\R_+$. The best solution to this (i.e. maximal $c$) is $c=1/4$ and $f(x)=\sqrt{x(1-x)}$. A more direct computation leads to a slightly stronger result.

\begin{corollary} \label{cor:chi-cheeger}
The spectral gap of a finite, irreducible Markov chain satisfies
$$
\~h \geq \lambda \geq 2\left(1-\sqrt{1-\~h^2/4}\right) \geq \frac{\~h^2}{4}\,.
$$
\end{corollary}

\begin{proof}
The upper bound is classical. The second lower bound follows from the first because $\sqrt{1-x}\leq 1-x/2$. For the first lower bound, consider $\C_{\sqrt{a(1-a)}}$ and apply equation \eqref{eqn:strong-Cheeger}. Then, for some $A\subset V$,
$$
\lambda \geq 
    2\left(1-
        \frac{\sqrt{1 + \~h(A)\pi(A^c)}\sqrt{1 - \~h(A)\pi(A)} 
            + \sqrt{1 - \~h(A)\pi(A^c)}\sqrt{1 + \~h(A)\pi(A)}} {2}
      \right)
$$
To simplify this let $X = \frac 12(1 + \~h(A)\,\pi(A^c))$
and $Y = \frac 12(1 - \~h(A)\,\pi(A))$ in Lemma \ref{lem:Appendix}.
\end{proof}

The Corollary is sharp on the two-point space $u-v$ with $\P(u,v)=\P(v,u)=1$ and $\~h=2$. Bounds in terms of $h$ typically show at best $\lambda\geq 1$ and so cannot be sharp for the two-point space.

Different choices of $f(a)$ work better if more is known about the dependence of $\~h(A)$ on set size. For example, for a walk on a cycle it is better to choose $f(a)=\sin(\pi a)$.

\begin{example} \label{ex:cycle}
Consider the reversible random walk on the cycle $C_n=\ZZ/n\ZZ$ of length $n$ with $\P(i,i\pm 1\mod n)=1/2$. If $A\subset C_n$ then $\Q(A,A^c)\geq 1/n$, and so the Cheeger inequality \eqref{eqn:normal-cheeger} gives the bound $\lambda\geq h^2/2=2/n^2$, while Corollary \ref{cor:chi-cheeger} improves this slightly to $\lambda\geq\~h^2/4=4/n^2$ or even the sharp $\lambda\geq 2$ when $n=2$. To apply Corollary \ref{cor:diffiQ} directly we solve the differential equation $-f/f''\leq c^{-1}$, or $f''+c\,f\geq 0$. The largest value of $c$ is obtained by the concave function $f(a)=\sin(\pi a)$ with $c=\pi^2$. Although $f''$ is not concave the function $f$
can still be used in Equation \eqref{eqn:strong-Cheeger} to obtain
$$
\lambda \geq 2\min_{\pi(A)\leq 1/2} 1-\cos(\pi\,\Q(A,A^c)) = 2(1-\cos(\pi/n)) \approx \frac{\pi^2}{n^2}\,.
$$
A more refined argument can be used to determine $\lambda$ exactly. As before, consider $\P'=\frac{I+\P}{2}$. Then
$$
\lambda_{\P} = 2\lambda_{\P'} \geq 1-\C_{\sin(\pi a),\P'}=1-\cos(2\pi/n)\,,
$$
which is the correct value of $\lambda$.
\end{example}

%********************** vertex-expansion ****************************

\subsection{Vertex-expansion} \label{sec:vertex}

The Generalized Cheeger inequality can also be used to show Cheeger-like inequalities in terms of vertex-expansion (the number of boundary vertices), leading to sharp versions of bounds due to Alon \cite{Alon86.1}, Bobkov, Houdr\'e and Tetali \cite{BHT00.1} and Stoyanov \cite{Sto01.1}.

Two notions of vertex-expansion are required:

\begin{definition}
If $A\subset V$ then the internal and external boundaries are $\partial_{in}(A)=\{x\in A:\,\Q(x,A^c)>0\}$ and $\partial_{out}(A)=\partial_{in}(A^c)=\{x\in A^c:\,\Q(x,A)>0\}$. The internal and external vertex expansion are
$$
h_{in}=\min_{\pi(A)\leq 1/2} h_{in}(A)
\qquad\textrm{and}\qquad
h_{out}=\min_{\pi(A)\leq 1/2} h_{out}(A)
$$
where
$$
h_{in}(A) = \frac{\pi(\partial_{in}(A))}{\pi(A)}
\qquad\textrm{and}\qquad
h_{out}(A) = \frac{\pi(\partial_{out}(A))}{\pi(A)}\,.
$$
Quantities $\~h_{in}$ and $\~h_{in}(A)$ are defined similarly, but with $\pi(A)\pi(A^c)$ in the denominator.
The minimum transition probability $\P_0 = \min_{x\neq y\in V} \{\P(x,y):\,\P(x,y)>0\}$ will also be required.
\end{definition}

\begin{theorem} \label{thm:vertex_congestion}
The spectral gap of a finite, reversible Markov kernel satisfies
\begin{eqnarray*}
\lambda &\geq&
           1-\sqrt{1-h_{out}\P_0} - \P_0\,\left(\sqrt{1+h_{out}}-1\right) 
  \geq \frac{\P_0}{12}\,\min\left\{h_{out}^2,\,h_{out}\right\} \\
\lambda &\geq&
           1-\sqrt{1+h_{in}\P_0}  - \P_0\,\left(\sqrt{1-h_{in}}-1\right) 
   \geq \frac{\P_0}{8}\,h_{in}^2 \\
\lambda &\geq& 1-\sqrt{1-\left(\frac{\~h_{in}\P_0}{2}\right)^2} 
             - \P_0\,\left(\sqrt{1-\left(\frac {\~h_{in}}{2}\right)^2}-1\right) 
   \geq \frac{\P_0(1+\P_0)}{8}\,\~h_{in}^2 
\,.
\end{eqnarray*}
For the non-reversible case replace $\P_0$ by $\P_0/2$. 
% $h_{in}$ by $h_{in}(\P)+h_{in}(\P^*)$ and $h_{out}$ by $h_{out}(\P)+h_{out}(\P^*)$.
\end{theorem}

\begin{proof}
First consider the reversible lazy case. Given $A\subset V$, a vertex $x\in A$ is in $\partial_{in}(A)$ if and only if $\P(x,A^c)>0$, which happens if and only if $\Q(A,x)\leq(1-\P_0)\pi(x)$, if and only if $x\notin A_u$ for every $u>1-\P_0$. Thus, given only $h_{in}(A)$, that $\pi(A_u)$ is non-increasing and that $\int_0^1\pi(A_u)\,du=\pi(A)$, then the integral $\int_0^t \pi(A_u)\,du$ is minimized for all $t\in[0,1]$ if $\pi(A_u)=m(u)$ is as in the first diagram of Figure \ref{fig:inside_boundary_extreme}.

\begin{figure}[ht]
\begin{center}
\includegraphics[height=1.25in]{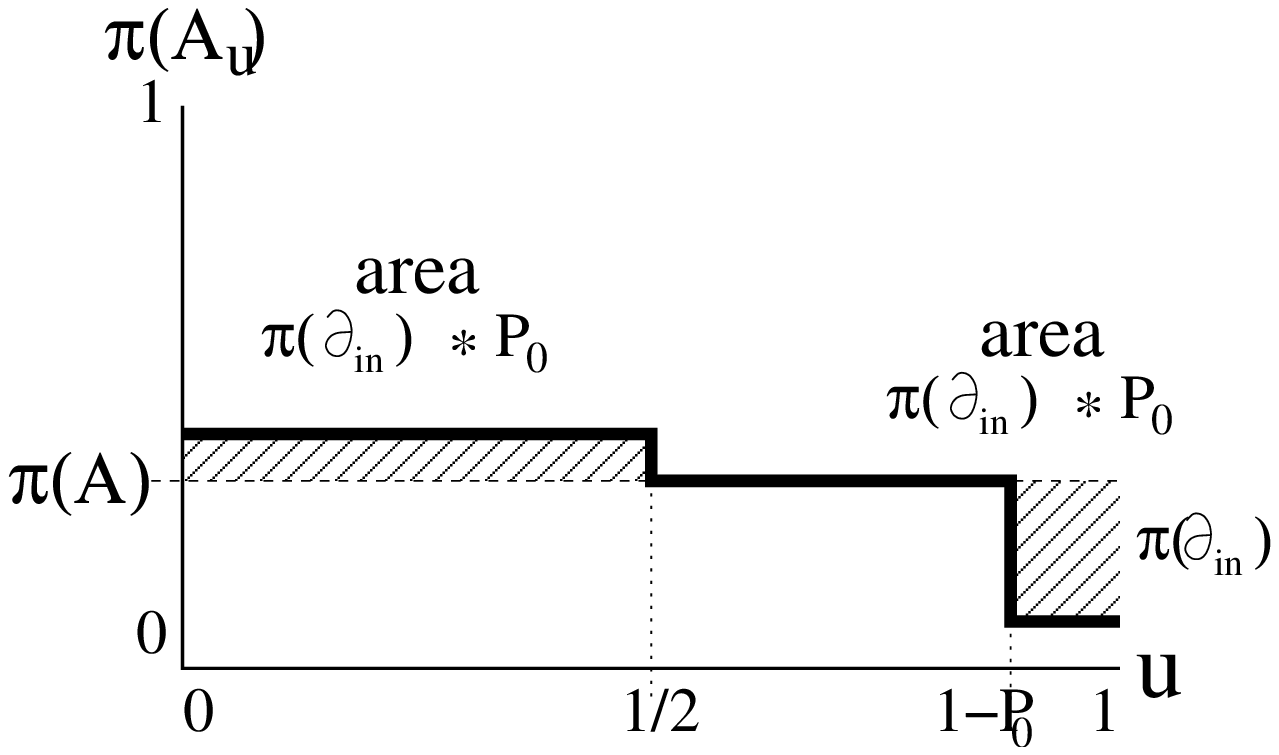}
\qquad
\includegraphics[height=1.25in]{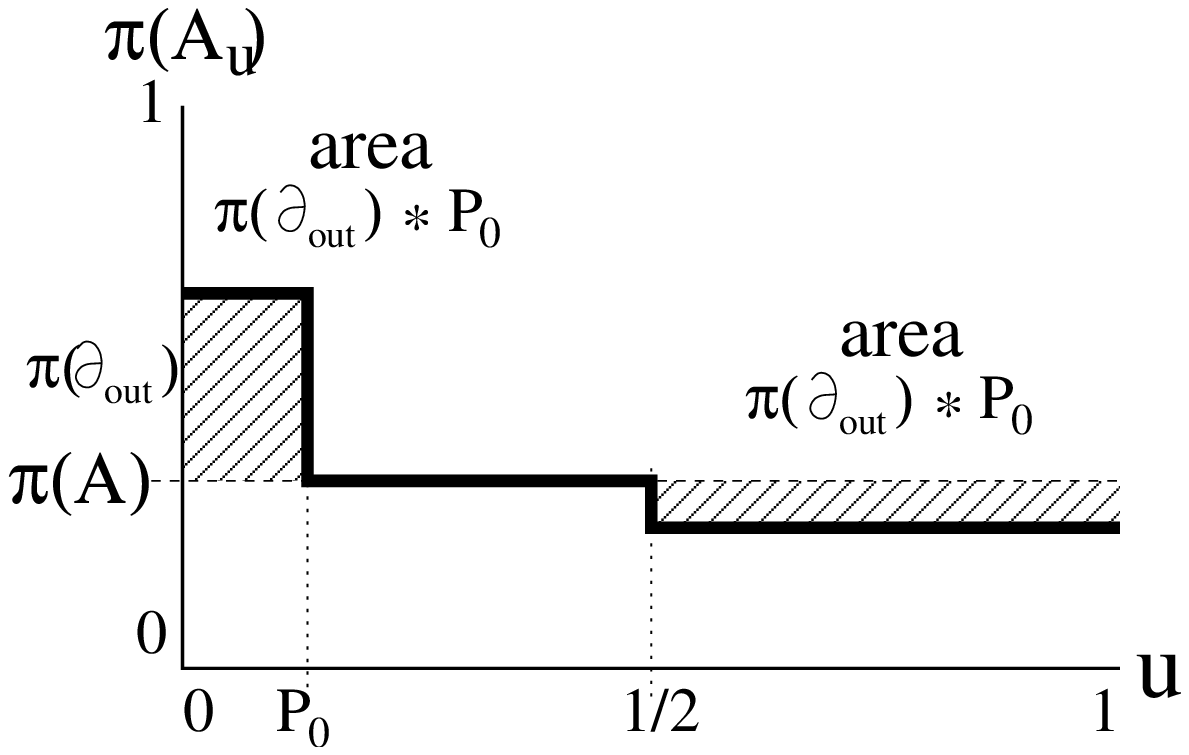}
\caption{Extreme case for $\P$ with $\partial_{in}(A)$, and $\P^*$ with $\partial_{out}(A)$.} \label{fig:inside_boundary_extreme}
\end{center}
\end{figure}

Finish by applying Lemma \ref{lem:worst-case} to upper bound $C_{\sqrt{a}}$ in terms of $h_{in}$ or $\C_{\sqrt{a(1-a)}}$ in terms of $\~h_{in}$. For $h_{out}$ upper bound $\C_{\sqrt{a}}$ using the second diagram of Figure \ref{fig:inside_boundary_extreme}. As before, the non-lazy case is reduced to the lazy case by the relation $\lambda_{\P}=2\lambda_{\P'}$ where $\P'=\frac{\I+\P}{2}$, while the non-reversible case is reduced to the reversible case by the relation $\lambda_{\P}=\lambda_{\P''}$ where $\P''=\frac{\P+\P^*}{2}$.
\end{proof}

To compare this to previous bounds we note that Stoyanov \cite{Sto01.1}, improving on results of Alon \cite{Alon86.1} and Bobkov, Houdr\'e and Tetali \cite{BHT00.1}, showed that a reversible Markov chain will satisfy
\begin{eqnarray*}
\lambda &\geq&
\max\left\{\frac{\P_0}{2}\,\left(1-\sqrt{1-h_{in}}\right)^2,\,
   \frac{\P_0}{4}\,\left(\sqrt{1+h_{out}}-1\right)^2\right\} \\
&\geq& \max\left\{\frac{\P_0}{8}\,h_{in}^2,\,\frac{\P_0}{24}\,\min\{h_{out}^2,h_{out}\}\right\}\,.
\end{eqnarray*}
Our Theorem \ref{thm:vertex_congestion}, and the approximations $\sqrt{1-h_{out}\P_0}\leq 1-h_{out}\P_0/2$ and $\sqrt{1+h_{in}\P_0}\leq 1+h_{in}\P_0/2$, give a stronger bound for reversible chains,
\begin{eqnarray*} \label{eqn:reversible_gap}
\lambda &\geq&
\frac{\P_0}{2}\max\left\{1-\sqrt{1-h_{in}},\,\sqrt{1+h_{out}}-1\right\}^2  \\
\lambda &\geq& \max\left\{ \frac{\P_0}{8}\,\~h_{in}^2,\,\frac{\P_0}{12}\min\{h_{out}^2,h_{out}\}\right\}\,.
\end{eqnarray*}

\begin{remark} \label{rmk:sharp}
The $h_{in}$ and $h_{out}$ bounds in this section were not sharp, despite our having promised sharp bounds. This is because $\C_{\sqrt{a(1-a)}}\leq\C_{\sqrt{a}}$ is a better quantity to consider. If $\C_{\sqrt{a(1-a)}}$ were used instead of $\C_{\sqrt{a}}$ then we would obtain sharp, although quite complicated, bounds; these bounds simplify in the $\tilde{h}$ and $\tilde{h}_{in}$ cases which is why we have used $\C_{\sqrt{a(1-a)}}$ for those two cases. Bounds based on $\C_{\sqrt{a(1-a)}}$ are sharp on the two-point space $u-v$ with $\P(u,v)=\P(v,u)=1$. 
\end{remark}

%********************* Mixing edge and vertex expansion *********************

\subsection{Mixing edge and vertex expansion}

We can easily combine edge and vertex-expansion quantities, and maximize at the set level rather than at a global level. For instance, in the reversible case
$$
\lambda \geq \min_{\pi(A)\leq 1/2} \max\left\{\frac 14\,\~h(A)^2,\, \frac{\P_0}{8}\,\~h_{in}(A)^2,\,\frac{\P_0}{12}\min\{h_{out}(A)^2,h_{out}(A)\}\right\}\,.
$$

Alternatively, we can apply Lemma \ref{lem:worst-case} directly:
\begin{theorem} \label{thm:various}
The spectral gap of a finite, reversible Markov kernel satisfies
\begin{eqnarray*}
\lambda \geq
 \min_{\pi(A)\leq 1/2} &&
              \,2 - \P_0\sqrt{1-h_{in}(A)} -\P_0\sqrt{1+h_{out}(A)}  \\
 &&- (1-\P_0)  \sqrt{1-\frac{h(A)-\P_0 h_{in}(A)}{1-\P_0}} 
   - (1-\P_0)\sqrt{1+\frac{h(A)-\P_0 h_{out}(A)}{1-\P_0}} \\
\end{eqnarray*}
$$
\lambda \geq \min_{\pi(A)\leq 1/2} 2-\P_0\sqrt{1-\left(\frac{\~h_{in}(A)}{2}\right)^2} 
                                          + \sqrt{1-\left(\frac{\~h(A)}{2}\right)^2} 
                + (1-\P_0)\,\sqrt{1 - \left(\frac{\~h(A)-\P_0\~h_{in}(A)}{2(1-\P_0)}\right)^2} \,.
$$
For the non-reversible case replace $\P_0$ by $\P_0/2$.
\end{theorem}

The proofs are no different from that of the cases already dealt with, other than that the worst cases $m(u)$ are somewhat more complicated, and so we omit the proofs. As in Remark \ref{rmk:sharp} the first bound can be made sharp (and even more complicated) by working with $1-\C_{\sqrt{a(1-a)}}$, while the second bound is already sharp on the two point space.

%\begin{proof}
%This is again similar to earlier cases, but this time the worst case for a lazy walk is just the function in Figure \ref{fig:mixed_extreme}, without the separate peak for $\partial_{out}^*$. The maximum on $\C_{\sqrt{a(1-a)}}$ occurs when $\pi(A)=1/2$, which gives a bound for a lazy chain, which generalizes by considering $\P'=\frac 12(I+\P)$, as before.
%\end{proof}

%******************************* Smallest eigenvalue ********************

\section{Bounding the smallest eigenvalue} \label{sec:small_eigenvalue}

The generalized Cheeger inequality can also be used to bound $1-\lambda_{max}$ for a reversible walk, by examining $\P$ directly instead of the lazy walk $\P'=\frac{\I+\P}{2}$ as before. Techniques of the previous sections carry through if modified expansion quantities are used, such as the following:

\begin{definition}
If $A\subset V$ then its {\em modified ergodic flow} is defined by
$$
\Psi(A) = \frac 12\,\int_0^1 |\pi(A_u)-\pi(A)|\,du\,.
$$
The {\em modified Cheeger constant} $\~\hbar$ is given by
$$
\~\hbar = \min_{A\subset V} \~\hbar(A)
\quad\textrm{where}\quad
\~\hbar(A) = \frac{\Psi(A)}{\pi(A)\pi(A^c)}\,.
$$
\end{definition}
By Lemma \ref{lem:ergodic-flow}, for a lazy chain $\Psi(A)=\Q(A,A^c)$ and hence also $\tilde{\hbar}(A)=\tilde{h}(A)$.

We can now show a lower bound on the eigenvalue gap:
\begin{theorem} \label{thm:profile}
Given a finite, irreducible Markov chain then
$$
1-\lambda_* \geq 1-\sqrt{1-\~\hbar^2} \geq \~\hbar^2/2\,.
$$
\end{theorem}

\begin{proof}
Let $\wp_A\in[0,1]$ be such that $\pi(A_u)\geq\pi(A)$ if $u<\wp_A$ and $\pi(A_u)\geq\pi(A)$ if $u>\wp_A$. Then $\Psi(A)=\int_0^{\wp_A} (\pi(A_u)-\pi(A))\,du=\int_{\wp_A}^1 (\pi(A)-\pi(A_u))\,du$ because $\int_0^1\pi(A_u)\,du=\pi(A)$. Apply Lemma \ref{lem:worst-case} with the second figure of Figure \ref{fig:extreme} to obtain
\begin{eqnarray*}
\C_{\sqrt{a(1-a)}}(A) &\leq& \sqrt{\left(\wp_A+\~\hbar(A)\,\pi(A^c)\right)    \left(\wp_A-\~\hbar(A)\,\pi(A)\right)} \\
     && + \sqrt{\left(1-\wp_A - \~\hbar(A)\,\pi(A^c)\right)\left(1-\wp_A + \~\hbar(A)\,\pi(A)\right)}
\end{eqnarray*}
To finish let $X = \displaystyle \wp_A + \~\hbar(A)\,\pi(A^c)$ and $Y = \displaystyle \wp_A - \~\hbar(A)\,\pi(A)$ in Lemma \ref{lem:Appendix}.
\end{proof}

For an isoperimetric interpretation of this, note that \cite{MT06.1} showed that
$$
\Psi(A) = \min_{\substack{B\subset \Omega,\,v\in \Omega,\\\pi(B)\leq\pi(A^c)<\pi(B\cup v)}}
              \Q(A,B)+\frac{\pi(A^c)-\pi(B)}{\pi(v)}\,\Q(A,v)\,.
$$
Hence, to bound spectral gap $\lambda$ consider the worst-case ergodic flow from a set $A$ to its complement $A^c$, whereas to bound $\lambda_*$ use the worst-case ergodic flow from a set $A$ to a set the same size as its complement $A^c$.

If we choose the $f$-congestion carefully then even better bounds may be achieved.

\begin{example}
Consider the cycle walk of Example \ref{ex:cycle}. For $x=\frac kn\leq \frac 12$ then $\min_{\pi(A)=x}\Psi(A)=\Q({\cal A}_x,{\cal B}_x)$ when ${\cal A}_x=\{0,2,4,\ldots,2k-2\}$ and ${\cal B}_x={\cal A}_x\cup\{-1,-2,-3,\ldots,-n+2k\}$. Then $\~\hbar=\Q({\cal A}_{1/2},{\cal B}_{1/2})=0$ if $n$ is even and so $1+\lambda_{n-1}\geq 0$, while $\~\hbar=\frac{2n}{n^2-1}\geq \frac 2n$ if $n$ is odd and so $1+\lambda_n\geq \~\hbar^2/2 \geq 2/n^2$. 

To improve on this, note that a bound similar to the lower bound of Corollary \ref{cor:diffiQ} holds for $\Psi(A)$ as well. Since $\Psi(A)\geq 1/2n$ for all $A\subset V$, this again suggests taking $f(a)=\sin(\pi a)$, and so if $n$ is odd then
$$
1-\lambda_{max} \geq 1-\C_{\sin(\pi a)} =
1-\cos(2\pi\Psi({\cal A}_{\frac{n-1}{2n}}))
=1-\cos \left(\frac{\pi}{n}\right)
$$
This is again an equality.
\end{example}

Vertex-expansion lower bounds for $1-\lambda_*$ (and hence also $1-\lambda_{max}$) hold as well. For instance, if $\hat\P_0=\min_{x,y\in V} \{\P(x,y):\,\P(x,y)>0\}$ (note that $x=y$ is permitted) and
$$
\hbar_{out} = \min_{\pi(A)\leq 1/2}\min_{\pi(B)=\pi(A^c)} \frac{\pi(\{x\in B:\,\Q(A,x)>0\})}{\pi(A)}
$$
then $1-\lambda_*\geq \frac{\hat\P_0}{12}\,\min\{\hbar_{out}^2,\,\hbar_{out}\}$. 

\begin{example}
A vertex-expander is a lazy walk where $h_{out}\geq\epsilon>0$. Analogously, we might define a non-lazy vertex-expander to be a walk where $\hbar_{out}\geq\epsilon>0$. If the expander is regular of degree $d$ then
$$
1-\lambda_{max} \geq \frac{\min\{\hbar_{out}^2,\,\hbar_{out}\}}{12d} \geq \frac{\epsilon^2}{12d}\,,
$$
which (up to a small constant factor) generalizes the relation $1-\lambda_{max}\geq\epsilon^2/4d$ for the lazy walk.
\end{example}

%*********************** Bibliography ********************************

\bibliographystyle{plain}
\bibliography{../references}

\begin{thebibliography}{1}

\bibitem{Alon86.1}
N.~Alon.
\newblock Eigenvalues and expanders.
\newblock {\em Combinatorica}, 6(2):83--96, 1986.

\bibitem{BHT00.1}
S.~Bobkov, C.~Houdr\'e, and P.~Tetali.
\newblock $\lambda_{\infty}$, vertex isoperimetry and concentration.
\newblock {\em Combinatorica}, 20(2):153--172, 2000.

\bibitem{DS91.1}
P.~Diaconis and D.~Stroock.
\newblock Geometric bounds for eigenvalues of markov chains.
\newblock {\em The Annals of Applied Probability}, 1:36--61, 1991.

\bibitem{JS88.1}
M.~Jerrum and A.~Sinclair.
\newblock Conductance and the rapid mixing property for markov chains: the
  approximation of the permanent resolved.
\newblock {\em Proceedings of the 20th Annual ACM Symposium on Theory of
  Computing (STOC 1988)}, pages 235--243, 1988.

\bibitem{LS88.1}
G.~Lawler and A.~Sokal.
\newblock Bounds on the $l^2$ spectrum for markov chains and markov processes:
  a generalization of cheeger's inequality.
\newblock {\em Transactions of the American Mathematical Society},
  309:557--580, 1988.

\bibitem{MT06.1}
R.~Montenegro and P.~Tetali.
\newblock {\em Mathematical Aspects of Mixing Times in Markov Chains}, volume
  1:3 of {\em Foundations and Trends in Theoretical Computer Science}.
\newblock NOW Publishers, Boston-Delft, June 2006.

\bibitem{MP05.1}
B.~Morris and Y.~Peres.
\newblock Evolving sets, mixing and heat kernel bounds.
\newblock {\em Probability Theory and Related Fields}, 133(2):245--266, 2005.

\bibitem{Sto01.1}
T.~Stoyanov.
\newblock {\em Isoperimetric and Related Constants for Graphs and Markov
  Chains}.
\newblock Ph.d.~thesis, Department of Mathematics, Georgia Institute of
  Technology, 2001.

\end{thebibliography}

%************************ Appendix *******************************

\section*{Appendix}

The following lemma was used for a few simplifications but was left for the Appendix.

\begin{lemma}\label{lem:Appendix}
If $X,Y\in[0,1]$ then
$$
\sqrt{XY}+\sqrt{(1-X)(1-Y)}\leq\sqrt{1-(X-Y)^2}\,.
$$ 
\end{lemma}

\begin{proof}
\begin{eqnarray*}
\left(\sqrt{X\,Y}+\sqrt{(1-X)(1-Y)}\right)^2 &=& 1-(X+Y)+2\,X\,Y \\
   && + \sqrt{[1-(X+Y)+2\,X\,Y]^2 - [1-2(X+Y)+(X+Y)^2]} \\
  &\leq& 2\left[1-(X+Y)+2\,X\,Y\right] - \left[1-2(X+Y)+(X+Y)^2\right] \\
       &=&    1+2\,X\,Y-X^2-Y^2 = 1-(X-Y)^2
\end{eqnarray*}
The inequality follows from the relation $\sqrt{a^2-b} \leq a-b$ if $a^2\geq b$, $a\leq\frac{1+b}{2}$ and $a\geq b$
(square both sides to show this), applied with $a=1-(X+Y)+2\,X\,Y$ and $b=1-2(X+Y)+(X+Y)^2$.
\end{proof}

\end{document}